\documentclass[a4paper,12pt]{amsart}
\usepackage{latexsym, amssymb,amsmath,amsthm}
\usepackage{mathrsfs}
\usepackage{bbm}
\usepackage[all, knot]{xy}
\xyoption{arc}

\usepackage{anysize}\marginsize{30mm}{30mm}{30mm}{30mm}

\def \1{\mathbf{1}}

\numberwithin{equation}{section}

\begin{document}

\date{}
\def\thefootnote{}

\title[Repetitive cluster-tilted algebras]{Repetitive cluster-tilted algebras$^\star$}

\author[S. Zhang]{Shunhua Zhang}
\address{School of Mathematics, Shandong University,
Jinan 250100, China} \email{shzhang@sdu.edu.cn}

\author[Y. Zhang]{Yuehui Zhang}
\address{Department of Mathematics,\ Shanghai JiaoTong University,\\
Shanghai 200240, P. R. China} \email{zyh@sjtu.edu.cn}
\date{}
\maketitle

\begin{abstract}
Let $H$ be a finite dimensional hereditary algebra over an
algebraically closed field $k$ and $\mathscr{C}_{F^m}$ be the
repetitive cluster category of $H$ with $m\geq 1$. We investigate
 the properties of cluster tilting objects in
$\mathscr{C}_{F^m}$ and the structure of repetitive cluster-tilted
algebras. Moreover, we generalized Theorem 4.2 in \cite{bmrrt} (Buan
A, Marsh R, Reiten I. Cluster-tilted algebra. Trans. Amer. Math.
Soc.,
 359(1)(2007), 323-332.) to the situation of $\mathscr{C}_{F^m}$,
 and prove that the tilting graph $\mathscr{K}_{\mathscr{C}_{F^m}}$
 of $\mathscr{C}_{F^m}$ is connected.
\end{abstract}

\footnote{ $^\star$ Supported by the NSFC(11171183).}

\vskip 0.3in

\section {Introduction}

Let $H$ be a finite dimensional hereditary algebra over an
algebraically closed field $k$. The endomorphism algebra of a
tilting module over $H$ is called tilted algebra. Cluster category
of type $H$ is the orbit category $\mathscr{C} =D^b(H)/(F)$ of the
derived category $D^b(H)$ of $H$ by an automorphism group generated
by $F =\tau^{-1}[1]$, where $\tau$ is the Auslander-Reiten
translation in $D^b(H)$ and $[1]$ is the shift functor of $D^b(H)$.
$\mathscr{C}$ is a triangulated category and is a Calabi-Yau
categories of CY-dimension $2$, see \cite{bmrrt, k}.  It was shown
that any cluster tilting object of $\mathscr{C}$ is induced by a
tilting module of a hereditary algebra $H'$ which is derived
equivalent to $H$. The endomorphism algebra of a tilting object in
$\mathscr{C}$ is called a cluster-tilted algebra.

\vskip 0.2in

Now cluster-tilted algebras and cluster categories provide an
algebraic understanding of combinatorics of cluster algebras defined
and studied by Fomin and Zelevinsky in \cite{fz}. In this
connection, the indecomposable exceptional objects in cluster
categories correspond to the cluster variables, and cluster tilting
objects (maximal 1-orthogonal subcategories \cite{i1, i2})
correspond to clusters in the corresponding cluster algebras, see
\cite{ck1, ck2}. Moreover by \cite{k} or \cite{kz}, cluster-tilted
algebras provide a class of Gorenstein algebras of Gorenstein
dimension $1$, which is important in representation theory of
algebras \cite{r1}.

\vskip 0.2in

For any positive integer $m$,  a repetitive cluster category
$\mathscr{C}_{F^m}=D^b(H)/(F^m)$, which is defined by \cite{z} as
the orbit category of the derived category $D^b(H)$ by the group
$(F^m)$ generated by $F^m$, is a triangulated category by Keller
\cite{k}, which is also a Calabi-Yau category of Calabi-Yau
dimension 2m/m. The cluster tilting objects in this repetitive
cluster category are shown to correspond one-to-one to those in the
classical cluster categories. The endomorphism algebras of cluster
tilting objects in $\mathscr{C}_{F^m}$ are called repetitive
cluster-tilted algebras. They all have the same representation type
and share a universal covering: the endomorphism algebra of
corresponding cluster tilting subcategory in $D^b(H)$, see \cite{z}
for details.

\vskip 0.2in

In this paper, we investigate the properties of cluster tilting
objects in $\mathscr{C}_{F^m}$ and the structure of repetitive
cluster-tilted algebras. The article is organized as follows: In
Section 2 we prove some basic facts about cluster tilting objects in
$\mathscr{C}_{F^m}$. In Section 3, we investigate the structure of
repetitive cluster-tilted algebras, and generalize Theorem 4.2 in
\cite{bmr} to the situation of $\mathscr{C}_{F^m}$, see Theorem 3.4.
Furthermore, we also prove that the tilting graph
$\mathscr{K}_{\mathscr{C}_{F^m}}$
 of $\mathscr{C}_{F^m}$ is connected.

\vskip 0.2in

Throughout this paper, we fix an algebraically closed field $k$, and
denote by  $H$ a finite dimensional hereditary $k$-algebra with $n$
simple modules. We denote by $\mathscr{C}=D^b(H)/(F)$ the cluster
category of $H$. A basic {\it tilting object} in $\mathscr{C}$ is an
object $T$ with $n$ non-isomorphic indecomposable direct summands
such that ${\rm Ext}^1_{\mathscr{C}}(T,T)=0$. We follow the standard
terminologies and notations used in the representation theory of
algebras, see \cite{ass, ars, h, r2}.

\vskip 0.3in

\section {Properties of cluster tilting objects in $\mathscr{C}_{F^m}$}

In this section, we first recall some definitions and collect some
known results which will be used later, then we prove some basic
facts about cluster tilting objects in $\mathscr{C}_{F^m}$.

\vskip 0.2in

Let $\mathscr{A}$ be a Krull-Remark-Schmidt category and
$\mathscr{B}$ a full subcategory of $\mathscr{A}$. We denote by
${\rm ind}\ \mathscr{B}$ the set of all indecomposable objects in
$\mathscr{B}$. For any object $M$ in $\mathscr{A}$, we denote by
${\rm add}\ M$ the full subcategory of $\mathscr{A}$ consisting of
finite direct sums of indecomposable summands of $M$ and by
$\delta(M)$ the number of non-isomorphic indecomposable summands of
$M$.

\vskip 0.2in

Let $m$ be an integer with $m\geq 1$ and $\mathscr{C}_{F^m}=
D^b(H)/(F^m)$ be the repetitive cluster category defined in [11]
with $F=\tau^{-1}[1]$.  Note that $\mathscr{C}_{F^1}= \mathscr{C}$
is the classical  cluster category. The following definition is
defined in \cite{i1, i2, z}.

\vskip 0.2in

{\bf Definition 2.1.} \  {\it An object $T$  of $\mathscr{C}_{F^m}$
 is a cluster tilting object provided $X\in {\rm add}\ T$ if and only
 if ${\rm Ext}^1_{\mathscr{C}_{F^m}}(X,T)=0$ and $X\in {\rm add}\ T$ if and only
 if ${\rm Ext}^1_{\mathscr{C}_{F^m}}(T,X)=0$.}

\vskip 0.2in

The triangle functor $\rho_m: \mathscr{C}_{F^m}\rightarrow
\mathscr{C}$ is defined in \cite{z} which is also a covering
functor. The following lemma is taken from \cite{z} which will be
used in the sequel.

\vskip 0.2in

{\bf Lemma 2.2.} \  (1). \ {\it \  $\mathscr{C}_{F^m}$ is a
Krull-remark-Schmidt category.}

\vskip 0.1in

 (2). \ {\it \  ${\rm ind}\ \mathscr{C}_{F^m}= \bigcup\limits_{i=0}^{m-1}
 ({\rm ind}\ F^i(\mathscr{C}))$.}

\vskip 0.1in

 (3). \ {\it \ $T$ is a cluster tilting object in $\mathscr{C}$ if
 and only if $\rho_m^{-1}(T)$ is a cluster tilting object in
$\mathscr{C}_{F^m}$.}

\vskip 0.1in

 (4). \ {\it \ For any tilting $H$-module $T$,
 $\bigoplus\limits_{i=0}^{m-1}F^iT$ is a cluster tilting object in
$\mathscr{C}_{F^m}$, and any cluster tilting object $M$ in
$\mathscr{C}_{F^m}$ arises in this way, i.e., there is a hereditary
algebra $H'$, which is derived equivalent to $H$, and a tilting
$H'$-module $T'$ such that the cluster tilting object $M$ is induced
from $T'$.}

 \vskip 0.2in

{\bf Remark.} \ It is easy to see that $\delta(T)=mn$ provided $T$
is a cluster tilting object in $\mathscr{C}_{F^m}$.

 \vskip 0.2in

Let $M$ be an object of $\mathscr{C}_{F^m}$. $M$ is said to be
$1$-orthogonal if ${\rm Ext}^1_{\mathscr{C}_{F^m}}(M,M)=0$. A basic
$1$-orthogonal object of $\mathscr{C}_{F^m}$ is said to be an almost
tilting object if $\delta(M)=nm-1$ and there exists an
indecomposable object $X$ of $\mathscr{C}_{F^m}$ such that $M\oplus
X$ is a cluster tilting object of $\mathscr{C}_{F^m}$.

\vskip 0.2in

{\bf Definition 2.3.} \  {\it An object $M$ of $\mathscr{C}_{F^m}$
is said to be $F$-stable if $M$ can be written as $M= X\oplus
FX\oplus\cdots \oplus F^{m-1}X$ for some object $X$ in
$\mathcal{C}$. In this case, we say that the $F$-stable object $M$
is determined by $X$.}

\vskip 0.2in

{\bf Lemma 2.4.} \  {\it Let $M$ be a $F$-stable object of
$\mathscr{C}_{F^m}$ and  $M= X\oplus FX\oplus\cdots \oplus F^{m-1}X$
for some object $X$ in $\mathcal{C}$. Then ${\rm
Ext}^1_{\mathscr{C}_{F^m}}(M,M)=0$ if and only if ${\rm
Ext}^1_{\mathscr{C}}(X,X)=0$.}

\vskip 0.1in

{\bf Proof.} \  \ It follows from that
\begin{eqnarray*}
{\rm Ext}^1_{\mathscr{C}_{F^m}}(M,M) &=& {\rm
Hom}_{\mathscr{C}_{F^m}}(M,M[1])\\
&=& \bigoplus\limits_{i\in \mathbb{Z}}{\rm Hom}_{\mathcal {D}^b(H)}(M,(F^m)^iM[1])\\
&=& \bigoplus\limits_{i\in \mathbb{Z}}{\rm Hom}_{\mathcal
{D}^b(H)}(\bigoplus\limits_{j=0}^{m-1}F^jX,F^{mi}
(\bigoplus\limits_{l=0}^{m-1}F^lX)[1])\\
&=&\bigoplus\limits_{m} \bigoplus\limits_{i\in \mathbb{Z}}{\rm Hom}_{\mathcal {D}^b(H)}(X,F^{mi}X[1]) \\
&=& \bigoplus\limits_{m}{\rm Hom}_{\mathscr{C}}(X,X[1])\\
&=&\bigoplus\limits_{m}{\rm Ext}^1_{\mathscr{C}}(X,X).
\end{eqnarray*}
Where $\bigoplus\limits_{m}X$ denotes the direct sum of $m$ copies
of $X$.      $\hfill\Box$

\vskip 0.2in

{\bf Lemma 2.5.}\ {\it  Let $X$ and $Y$ be objects in $\mathscr{C}$,
and $M=X\oplus FX\oplus\cdots \oplus F^{m-1}X$ be the object in
$\mathscr{C}_{F^m}$ determined by $X$. Then ${\rm
Hom}_{\mathscr{C}_{F^m}}(M,F^jY)\simeq {\rm
Hom}_{\mathscr{C}_{F^m}}(M,Y)$ for $0\leq j \leq m-1$.}

\vskip 0.1in

{\bf Proof.}\  It follows from that $F^jM \simeq M$ in
$\mathscr{C}_{F^m}$.           $\hfill\Box$

\vskip 0.2in

{\bf Proposition 2.6.}  \ {\it Let $M$ be a basic almost tilting
object in $\mathscr{C}_{F^m}$ with $m\geq 2$. Then $M$ has only one
indecomposable complement $X$ in $\mathscr{C}_{F^m}$.}

\vskip 0.1in

{\bf Proof.}\  Let $X$ be an indecomposable complement to $M$.
Assume that $Y$ is another indecomposable complement to $M$, we want
to show that $Y\simeq X$. Since $M\oplus X$ is a basic tilting
object in $\mathscr{C}_{F^m}$, we may assume that $M\oplus X$ is
determined by a basic tilting object $T$ in $\mathscr{C}$, that is,
$M\oplus X= T\oplus FT\oplus \cdots \oplus F^{m-1}T$. Note that
$m\geq 2$ and $\delta(M)=nm-1$, we have that $T\oplus FT\oplus
\cdots \oplus F^{m-1}T=M\oplus Y$. Then $Y\simeq X$ follows from
Lemma 2.2 since $\mathscr{C}_{F^m}$ is a Krull-Remark-Schmidt
category.      $\hfill\Box$

\vskip 0.2in

Let $M$ be the $F$-stable object in $\mathscr{C}_{F^m}$ determined
by an object $X$ in $\mathscr{C}$. We denote by $O(M)$ the
 number of $F$-orbit in $M$ determined by the indecomposable summands of
 $M$. It is easy to see that $O(M)=\delta (M)$.

\vskip 0.2in

{\bf Lemma 2.7.} \  {\it Let $M$ be a $F$-stable object of
$\mathscr{C}_{F^m}$ and  ${\rm Ext}^1_{\mathscr{C}_{F^m}}(M,M)=0$.
Then $M$ is a cluster tilting object of $\mathscr{C}_{F^m}$ if and
only if $O(M)=n$.}

\vskip 0.1in

{\bf Proof.} \  Assume that $M$ is determined by $X$ in
$\mathscr{C}$, that is, $M=X\oplus FX\oplus\cdots \oplus F^{m-1}$.
According to Lemma 2.4, we have that ${\rm
Ext}^1_{\mathscr{C}}(X,X)=0$. It is well known that $X$ is a tilting
object in $\mathscr{C}$ if and only if $\delta(X)=n$.  The
consequence follows from Lemma 2.2.         $\hfill\Box$

\vskip 0.2in

{\bf Definition 2.8.} \  {\it An object $M$ of $\mathscr{C}_{F^m}$
is said to be rigid if $M$ is $F$-stable and ${\rm
Ext}^1_{\mathscr{C}_{F^m}}(M,M)=0$.  A rigid object $M$ of
$\mathscr{C}_{F^m}$ is said to be an almost near tilting object if
$O(M)=n-1$.}

\vskip 0.2in

{\bf Proposition 2.9.} \  {\it Let $M$ be an almost near tilting
object of $\mathscr{C}_{F^m}$. Then $M$ has exactly two kinds of
complements. That is, there exist two $F$-stable objects $M_1$ and
$M_2$, determined by non-isomorphic indecomposable objects $X_1$ and
$X_2$ of $\mathscr{C}$ respectively, such that $M\oplus M_1$ and
$M\oplus M_2$ are cluster tilting objects in $\mathscr{C}_{F^m}$.}

\vskip 0.1in

{\bf Proof.} \  Assume that $M$ is determined by an object $X$ in
$\mathscr{C}$ and  $M=X\oplus FX\oplus\cdots \oplus F^{m-1}X$. Then
$\delta (X)=n-1$, and by using Lemma 2.4 we have that ${\rm
Ext}^1_{\mathscr{C}}(X,X)=0$. Therefore, $X$ is an almost tilting
object of $\mathscr{C}$. According to \cite{bmrrt}, $X$ has exactly
two non-isomorphic complements $X_1$ and $X_2$ in $\mathscr{C}$.
Assume that $M_i=X_i\oplus FX_i\oplus\cdots \oplus F^{m-1}X_i$ for
$i=1,2$, then it is easy to see that $M\oplus M_1$ and $M\oplus M_2$
are cluster tilting objects in $\mathscr{C}_{F^m}$. $\hfill\Box$

\vskip 0.3in

\section {Repetitive cluster-tilted algebras}

Let $M$ be a cluster tilting object in $\mathscr{C}_{F^m}$. Then the
endomorphism algebra ${\rm End}_{\mathscr{C}_{F^m}}(M)$ is called a
repetitive cluster-tilted algebra.

\vskip 0.2in

{\bf Proposition 3.1.} \  {\it Let $T$ be a basic tilting module of
$H$ and $M$ be the $F$-stable object in $\mathscr{C}_{F^m}$
determined by $T$. Then $M$ is a cluster tilting object in
$\mathscr{C}_{F^m}$ and the endomorphism algebra ${\rm
End}_{\mathscr{C}_{F^m}}(M)$ is isomorphic to
$$\begin{pmatrix}
                  C_{0} &  &  &   \\
                  E_{1} & C_{1} &   \\
                   & \ddots & \ddots &  \\
                  & & E_{m-1} & C_{m-1} \\
                \end{pmatrix}.
$$
Where $C_{i}=C= {\rm End}_H\ T$ and $E_{i}= {\rm Ext}^2_C(DC,C)$ for
$0\leq i \leq m-1$, all the remaining coefficients are zero and
multiplication is induced from the canonical isomorphisms
$C\otimes_C E\cong \ _CE\cong E\otimes_CC$  and the zero morphism
$E\otimes_C E \longrightarrow 0$.}

\vskip 0.1in

{\bf Proof.} \ By the assumption, $M=T\oplus FT\oplus\cdots \oplus
F^{m-1}T$.  As a vector space, we have
$$
{\rm End}_{\mathscr{C}_{F^m}}(M)= \bigoplus\limits_{i\in
\mathbb{Z}}{\rm Hom}_{\mathcal
{D}^b(H)}(\bigoplus\limits_{j=0}^{m-1}F^jT,F^{mi}(\bigoplus\limits_{l=0}^{m-1}F^lT)).
$$
Since $T$ is a $H$-module, we have that ${\rm Hom}_{\mathcal
{D}^b(H)}(F^iT,F^jT)= 0$ unless $i=j$ or $i=j-1$.  Moreover, ${\rm
Hom}_{\mathcal {D}^b(H)}(F^iT,F^iT)= {\rm Hom}_H (T,T)=C$ and ${\rm
Hom}_{\mathcal {D}^b(H)}(F^iT,F^{i+1}T)= {\rm Hom}_{\mathcal
{D}^b(H)} (T,FT)= {\rm Ext}^2_C(DC, C)$.         $\hfill\Box$

\vskip 0.2in

{\bf Remark.}\ According to Lemma 2.2.(4),  every repetitive
cluster-tilted algebra can be described as in Proposition 3.1.

\vskip 0.2in

{\bf Lemma 3.2.} \ {\it  Let $M$ be an indecomposable object in
$\mathscr{C}_{F^m}$ with  $m\geq 2$. Assume that $M=F^jX$ with $X$,
an indecomposable object in $\mathscr{C}$,  and that  $0\leq j\leq
m-1$. If $\widehat{M}= F^j X\oplus \cdots\oplus F^{m-1}X \oplus
F^{j-1}X\oplus \cdots\oplus X$ is a rigid object in
$\mathscr{C}_{F^m}$, then ${\rm End}_{\mathscr{C}_{F^m}}(M)$ is a
field.}

\vskip 0.1in

{\bf Proof.} \ By assumption, we have that ${\rm
Ext}^1_{\mathscr{C}_{F^m}}(\widehat{M}, \widehat{M})=0$. According
to Lemma 2.4, we have that ${\rm Ext}^1_{\mathscr{C}}(X, X)=0$.
Hence ${\rm End}_{\mathcal {D}^b(H)}(X,X) \simeq k$, since $m\geq
2$. We have that
\begin{eqnarray*}
{\rm End}_{\mathscr{C}_{F^m}}(M) &=& {\rm
End}_{\mathscr{C}_{F^m}}(X)\\
&=& \bigoplus\limits_{i\in \mathbb{Z}}{\rm Hom}_{\mathcal {D}^b(H)}(X,F^{mi}X)\\
&=& {\rm Hom}_{\mathcal{D}^b(H)}(X,X) =  k.
\end{eqnarray*}
  $\hfill\Box$

\vskip 0.2in

Let $M'$ be a basic almost near tilting object of
$\mathscr{C}_{F^m}$. We may assume that $M'=X\oplus FX\oplus\cdots
\oplus F^{m-1}X$ with $X$ be a basic almost tilting object in
$\mathscr{C}$. Let $N_1$ and $N_2$ be non-isomorphic $F$-stable
complements of $M'$. Then $M_1=M'\oplus N_1$ and $M_2=M'\oplus N_2$
are cluster tilting objects in $\mathscr{C}_{F^m}$. We denote by
$\Lambda_i={\rm End}_{\mathscr{C}_{F^m}}(M_i)$ for $i=1,2$. If $N_1$
and $N_2$ are determined by non-isomorphic indecomposable objects
$X_1$ and $X_2$ of $\mathscr{C}$ respectively, then $N_1= X_1\oplus
FX_1\oplus\cdots \oplus F^{m-1}X_1$ and $N_2= X_2\oplus
FX_2\oplus\cdots \oplus F^{m-1}X_2$.

\vskip 0.2in

{\bf Lemma 3.3.} \ {\it  Take the notation as above. Then
$S_{12}={\rm Hom}_{\mathscr{C}_{F^m}}(M_1, N_2[1])$ is a  semisimple
$\Lambda_1$-module and $S_{21}={\rm Hom}_{\mathscr{C}_{F^m}}(M_2,
N_1[1])$ is a semisimple  $\Lambda_2$-module.}

\vskip 0.1in

{\bf Proof.} \ By duality, we only need to prove that $S_{12}$ is a
 semisimple $\Lambda_1$-module.

For an integer $j$ with $0\leq j\leq m-1$, by Lemma 2.5 we have the
following isomorphism of $k$-spaces:  $ {\rm
Hom}_{\mathscr{C}_{F^m}}(M_1, F^jX_2[1])\simeq {\rm
Hom}_{\mathscr{C}_{F^m}}(M_1, X_2[1])$.

We claim that ${\rm Hom}_{\mathscr{C}_{F^m}}(M_1, X_2[1])$ is a
simple $\Lambda_1$-module.

In fact, if $m=1$, our claim is proved in \cite{bmr}.

If $m\geq 2$, we consider the following triangle $X_2 \rightarrow  B
\rightarrow X_1  \rightarrow X_2[1]$, where $B \rightarrow X_1 $ is
the minimal right ${\rm add}\ X$-approximation. Applying ${\rm
Hom}_{\mathscr{C}_{F^m}}(M_1, -)$ we obtain the following exact
sequence ${\rm Hom}_{\mathscr{C}_{F^m}}(M_1, B) \rightarrow {\rm
Hom}_{\mathscr{C}_{F^m}}(M_1, X_1) \rightarrow {\rm
Hom}_{\mathscr{C}_{F^m}}(M_1, X_2[1]) \rightarrow 0$

\begin{eqnarray*}
{\rm Hom}_{\mathscr{C}_{F^m}}(M_1, X_2[1]) &=& {\rm
Hom}_{\mathscr{C}_{F^m}}(M'\oplus N_1, X_2[1])\\
&=& {\rm Hom}_{\mathscr{C}_{F^m}}(N_1, X_2[1]) \\
&=& {\rm Hom}_{\mathscr{C}_{F^m}}(X_1\oplus
FX_1\oplus\cdots \oplus F^{m-1}X_1, X_2[1]) \\
&=& \bigoplus\limits_{j\in \mathbb{Z}}{\rm
Hom}_{\mathcal{D}^b(H)}(X_1\oplus
FX_1\oplus\cdots \oplus F^{m-1}X_1, F^{jm}X_2[1]) \\
&=& {\rm Hom}_{\mathcal{D}^b(H)}(X_1,X_2[1])\oplus {\rm
Hom}_{\mathcal{D}^b(H)}(X_1, FX_2[1])\\
&=& {\rm Hom}_{\mathscr{C}}(X_1,X_2[1]).
\end{eqnarray*}

Since $k$ is an algebraically closed field, by \cite{bmrrt}, ${\rm
Hom}_{\mathscr{C}}(X_1,X_2[1])$ is one-dimensional $k$-space. Thus
${\rm Hom}_{\mathscr{C}_{F^m}}(M_1, X_2[1])$ is a simple
$\Lambda_1$-module, our claim is proved.

By Lemma 2.5 again, we have that $S_{12}= {\rm
Hom}_{\mathscr{C}_{F^m}}(M_1, N_2[1])$ is a semisimple
$\Lambda_1$-module.  This completes the proof.      $\hfill\Box$

\vskip 0.2in

According to Corollary 4.4 in \cite{kz}, we know that ${\rm mod}\
\Lambda_1 \simeq \mathscr{C}_{F^m}/ {\rm add}\ M_1[1]$ and ${\rm
mod}\ \Lambda_2 \simeq \mathscr{C}_{F^m}/ {\rm add}\ M_2[1]$. Let
$\widetilde{M}= M'\oplus N_1\oplus N_2$.  We denote by $S_{12}$ the
semisimple $\Lambda_1$-module ${\rm Hom}_{\mathscr{C}_{F^m}}(M_1,
N_2[1])$ and by $S_{21}$ the semisimple $\Lambda_2$-module ${\rm
Hom}_{\mathscr{C}_{F^m}}(M_2, N_1[1])$, then we get equivalences
${\rm mod}\ \Lambda_1/{\rm add}\ S_{12} \simeq \mathscr{C}_{F^m}/
{\rm add}\ \widetilde{M}[1]$ and ${\rm mod}\ \Lambda_2/{\rm add}\
S_{21} \simeq \mathscr{C}_{F^m}/ {\rm add}\ \widetilde{M}[1]$.

\vskip 0.2in

Summarizing the above discussions, we get the following theorem
which is a generalization of Theorem 4.2 in \cite{bmr}.

\vskip 0.2in

{\bf Theorem 3.4.} \ {\it Let $M'$ be a basic almost near tilting
object of $\mathscr{C}_{F^m}$ determined by a basic almost tilting
object $X$ in $\mathscr{C}$ with non-isomorphic complements $X_1$
and $X_2$ be of  in $\mathscr{C}$. Then  $N_1= X_1\oplus
FX_1\oplus\cdots \oplus F^{m-1}X_1$ and $N_2= X_2\oplus
FX_2\oplus\cdots \oplus F^{m-1}X_2$ are non-isomorphic $F$-stable
complements of $M'$. Moreover, $M_1=M'\oplus N_1$ and $M_2=M'\oplus
N_2$ are cluster tilting objects in $\mathscr{C}_{F^m}$. Let
$\Lambda_1={\rm End}_{\mathscr{C}_{F^m}}(M_1)$ and $\Lambda_2={\rm
End}_{\mathscr{C}_{F^m}}(M_2)$. Then ${\rm mod}\ \Lambda_1/{\rm
add}\ S_{12} \simeq {\rm mod}\ \Lambda_2/{\rm add}\ S_{21}$, where
$S_{12}={\rm Hom}_{\mathscr{C}_{F^m}}(M_1, N_2[1])$ and $S_{21}={\rm
Hom}_{\mathscr{C}_{F^m}}(M_2, N_1[1])$ are  semisimple modules.}

\vskip 0.2in

Let $\mathcal {T}_{\mathscr{C}_{F^m}}$ be the set of all basic
cluster tilting objects in $\mathscr{C}_{F^m}$ up to isomorphism.
According to \cite{hu}, the tilting graph
$\mathscr{K}_{\mathscr{C}_{F^m}}$ of $\mathscr{C}_{F^m}$ is defined
as the following. The vertices of $\mathscr{K}_{\mathscr{C}_{F^m}}$
are the elements of $\mathcal {T}_{\mathscr{C}_{F^m}}$. There is an
edge between $M'$ and $M$ if there exists an almost near tilting
object $B$ such that $M'=B\oplus N_1$ and $M=B\oplus N_2$ with $N_1$
and $N_2$ are determined by indecomposables $X$ and $Y$ in
$\mathscr{C}$. That is, $B$ is determined by an almost tilting
object $T$ in $\mathscr{C}$ and $T$ has exactly two non-isomorphic
indecomposable complements $X$ and $Y$, such that $N_1$ and $N_2$
are determined by $X$ and $Y$ respectively.

\vskip 0.2in

 {\bf Theorem 3.5.} {\it The tilting graph
 $\mathscr{K}_{\mathscr{C}_{F^m}}$  of  $\mathscr{C}_{F^m}$ is connected.}

\vskip 0.1in

{\bf Proof.} Let $M_1$ and $M_2$ be two elements in $\mathcal
{T}_{\mathscr{C}_{F^m}}$. We suppose that  $M_1$ and $M_2$ are
determined  by basic tilting objects $T_1$ and $T_2$ of
$\mathscr{C}$ respectively. According to \cite{bmrrt} the tilting
graph $\mathscr{K}_{\mathscr{C}}$ of $\mathscr{C}$ is connected,
hence there exist basic tilting objects $X_1,\cdots,X_t$ of
$\mathscr{C}$ such that there is a path $T_1-X_1-\cdots-X_t-T_2$ in
tilting graph $\mathscr{K}_{\mathscr{C}}$ of $\mathscr{C}$. We
denote by $N_i$ the element of $\mathcal {T}_{\mathscr{C}_{F^m}}$
determined by $X_i$, i.e., $N_i= X_i\oplus FX_i\oplus\cdots \oplus
F_{m-1}X_i$ for $1\leq i\leq t$, according to Proposition 2.9,  we
obtain a path $M_1-N_1-\cdots-N_t-M_2$ in tilting graph
$\mathscr{K}_{\mathscr{C}_{F^m}}$. The proof is completed.
$\hfill\Box$

\vskip 0.2in

{\bf Acknowledgments.}\ \  We would like to thank the referee for
his or her valuable comments and suggestions, which improve the
presentation of this paper.

\end{document}